\documentclass[12pt,letterpaper,reqno]{amsart}

\usepackage{times}
\usepackage[T1]{fontenc}
\usepackage{mathrsfs}
\usepackage{latexsym}
\usepackage[dvips]{graphics}
\usepackage{epsfig}
\usepackage{amsmath,amsfonts,amsthm,amssymb,amscd}
\input amssym.def
\input amssym.tex

\newcommand\be{\begin{equation}}
\newcommand\ee{\end{equation}}
\newcommand\bea{\begin{eqnarray}}
\newcommand\eea{\end{eqnarray}}
\newcommand\bi{\begin{itemize}}
\newcommand\ei{\end{itemize}}
\newcommand\ben{\begin{enumerate}}
\newcommand\een{\end{enumerate}}
\newcommand\bc{\begin{center}}
\newcommand\ec{\end{center}}
\newcommand\ba{\begin{array}}
\newcommand\ea{\end{array}}

\newcommand{\N}{\mathbb{N}}
\newcommand{\E}{\mathbb{E}}

\newcommand{\foh}{\frac{1}{2}}  

\newtheorem{thm}{Theorem}[section]
\newtheorem{conj}[thm]{Conjecture}

\newtheorem{lem}[thm]{Lemma}

\newtheorem{defi}[thm]{Definition}

\theoremstyle{definition}
\newtheorem{rek}[thm]{Remark}

\newcommand{\ncr}[2]{{#1 \choose #2}}
\newcommand{\twocase}[5]{#1 \begin{cases} #2 & \text{{\rm #3}}\\ #4
&\text{{\rm #5}} \end{cases}   }

\begin{document}

\title{When almost all sets are
difference dominated}

\author{Peter Hegarty}
\email{hegarty@math.chalmers.se} \address{Mathematical Sciences,
Chalmers University Of Technology and G\"oteborg University,
G\"oteborg, Sweden}

\author{Steven J. Miller}
\email{Steven.J.Miller@williams.edu} \address{Department of Mathematics and
Statistics,
Williams College, Williamstown, MA 01267}

\subjclass[2000]{11P99 (primary), 11K99 (secondary).} \keywords{Sum
dominated sets, Binary linear forms, Strong concentration, Thresholds.}

\date{\today}


\begin{abstract} We investigate the relationship between the sizes of the
sum and difference sets attached to a subset of $\{0,1,...,N\}$,
chosen randomly according to a binomial model with parameter $p(N)$,
with $N^{-1} = o(p(N))$. We show that the random subset is almost
surely difference dominated, as $N \rightarrow \infty$, for any
choice of $p(N)$ tending to zero, thus confirming a conjecture of
Martin and O'Bryant. The proofs use recent strong concentration
results.
\par Furthermore, we exhibit a threshold phenomenon regarding the
ratio of the size of the difference- to the sumset. If $p(N) =
o(N^{-1/2})$ then almost all sums and differences in the random
subset are almost surely distinct, and in particular the difference
set is almost surely about twice as large as the sumset. If
$N^{-1/2} = o(p(N))$ then both the sum and difference sets almost
surely have size $(2N+1) - O(p(N)^{-2})$, and so the ratio in
question is almost surely very close to one. If $p(N) = c \cdot
N^{-1/2}$ then as $c$ increases from zero to infinity (i.e., as the
threshold is crossed), the same ratio almost surely decreases
continuously from two to one according to an explicitly given
function of $c$.
\par We also extend our results to the comparison of
the generalized difference sets attached to an arbitrary pair of
binary linear forms. For certain pairs of forms $f$ and $g$, we show
that there in fact exists a sharp threshold at $c_{f,g} \cdot
N^{-1/2}$, for some computable constant $c_{f,g}$, such that one
form almost surely dominates below the threshold, and the other
almost surely above it.
\par The heart of our approach involves using different tools to
obtain strong concentration of the sizes of the sum and difference sets
about their mean values, for various ranges of the parameter $p$.
\end{abstract}


\maketitle

\setcounter{equation}{0}

\setcounter{equation}{0}

\section{Introduction}

To know whether a random variable is strongly concentrated is an issue of
fundamental importance in many areas of mathematics and statistics.
In this paper we apply recent
results of Kim and Vu \cite{KiVu,Vu1,Vu2} to completely solve a combinatorial
number theory
question on the size of difference- and sumsets of integers. A
classical strong concentration
result (due to Chernoff) states that if $Y = \sum_{i=1}^n
t_i$ with the $t_i$ i.i.d. binary random variables, then for any
$\lambda > 0$ we have ${\rm Prob}(|Y - \E[Y]| \ge \sqrt{\lambda n})
\le 2e^{-\lambda/2}$. Within number theory, this result was used by
Erd\H{o}s (see \cite{AS}, Chapter 8) to prove the existence of so-called
$\lq$thin' bases of $\mathbb{N}$ of order 2. The general requirement
for many applications is to obtain Chernoff-like
exponential deviation bounds in situations when the atom variables
$t_{i}$ are not independent. For modern surveys of strong concentration
inequalities see, for example,
\cite{Ta} and \cite{Vu2}; the latter, in particular, contains a fine
selection of applications in random graph theory,
combinatorial number theory and finite geometry.

The specific result we shall utilise is a martingale inequality which
appears as Lemma 3.1 in \cite{Vu2}. It is an extension of the
classical Azuma inequality (\cite{AS}, Chapter 7) to functions whose
Lipschitz coefficients are small $\lq$on average'.
As remarked in \cite{Vu2}, this type of inequality is very general
and robust, and is expected to be applicable in numerous
situations; this is definitely true for our problem.

Let $S$ be a
subset of the integers. We define the sumset $S+S$ and difference
set $S-S$ by \bea S+S & \ = \ & \{s_1+s_2: s_i \in S\} \nonumber\\
S-S &=& \{s_1 - s_2: s_i \in S\}, \eea and denote the cardinality of
a set $A$ by $|A|$. As addition is commutative and subtraction is
not, a typical pair of integers generates two differences but only
one sum. It is therefore reasonable to expect a generic finite set
$S$ will have a larger difference set than sumset. We say a set is
\emph{sum dominated} (such sets are also called \emph{more sums than
differences}, or MSTD, sets) if the cardinality of its sumset
exceeds that of its difference set. If the two cardinalities are
equal we say the set is \emph{balanced}, otherwise \emph{difference
dominated}. Sum dominated sets exist: consider for example
$\{0,2,3,4,7,11,12,14\}$ (see \cite{He,MS,Na2} for additional
examples). In \cite{Na1}, Nathanson wrote \emph{``Even though there
exist sets $A$ that have more sums than differences, such sets
should be rare, and it must be true with the right way of counting
that the vast majority of sets satisfies $|A-A| > |A+A|$.''}

Recently Martin and O'Bryant \cite{MO} showed there are many sum
dominated sets. Specifically, let $I_N = \{0,\dots,N\}$. They prove
the existence of a universal constant $\kappa_{{\rm SD}} > 0$ such
that, for any $N \geq 14$, at least $\kappa_{{\rm SD}} \cdot 2^{N+1}$
subsets of $I_{N}$ are sum dominated (there are no sum dominated sets in
$I_{13}$). Their proof is based on
choosing a subset of $I_N$ by picking each $n \in I_N$ independently
with probability $1/2$. The argument can be generalized to
independently picking each $n\in I_N$ with any probability $p \in
(0,1)$, and yields the existence of a constant $\kappa_{{\rm SD},p}
> 0$ such that, as $N\to\infty$, a randomly chosen (with

respect to this model) subset is sum dominated with probability at
least $\kappa_{{\rm SD},p}$. Similarly one can prove there are
positive constants $\kappa_{{\rm DD},p}$ and $\kappa_{{\rm B},p}$
for the probability of having a difference dominated or balanced
set.

While the authors remark that, perhaps contrary to intuition, sum
dominated sets are ubiquitous, their result is a consequence of how
they choose a probability distribution on the space of subsets of
$I_N$. Suppose $p = 1/2$, as in their paper. With high probability a
randomly chosen subset will have $N/2$ elements (with errors of size
$\sqrt{N}$). Thus the density of a generic subset to the underlying
set $I_N$ is quite high, typically about $1/2$. Because it is so
high, when we look at the sumset (resp., difference set) of a
typical $A$ there are many ways of expressing elements as a sum
(resp., difference) of two elements of $A$. For example (see
\cite{MO}), if $k \in A+A$ then there are roughly $N/4 - |N-k|/4$
ways of writing $k$ as a sum of two elements in $A$ (similarly, if
$k\in A-A$ there are roughly $N/4-|k|/4$ ways of writing $k$ as a
difference of two elements of $A$). This enormous redundancy means
almost all numbers which can be in the sumset or difference set are.
In fact, using uniform density on the subsets of $I_N$ (i.e., taking
$p = 1/2$), Martin and O'Bryant show that the average value of
$|A+A|$ is $2N-9$ and that of $|A-A|$ is $2N-5$ (note each set has
at most $2N+1$ elements). In particular, it is only for $k$ near
extremes that we have high probability of not having $k$ in an $A+A$
or an $A-A$. In \cite{MO} they prove a positive percentage of
subsets of $I_N$ (with respect to the uniform density) are sum
dominated sets by specifying the fringe elements of $A$. Similar
conclusions apply for any value of $p > 0$.

At the end of their paper, Martin and O'Bryant conjecture that if, on the
other hand, the parameter $p$ is a function of $N$ tending to zero
arbitrarily slowly, then
as $N \rightarrow \infty$ the probability that a randomly chosen subset of
$I_{N}$ is sum dominated should also tend to zero. In this paper we will,
among other things, prove this conjecture.

We shall find it convenient to adopt the following (fairly standard)
shorthand notations. Let $X$ be a real-valued random variable
depending on some positive integer parameter $N$, and let $f(N)$ be
some real-valued function. We write $\lq X \sim f(N)$' to denote the
fact that, for any $\epsilon_{1}, \epsilon_{2} > 0$, there exists
$N_{\epsilon_{1},\epsilon_{2}} > 0$ such that, for all $N >
N_{\epsilon_{1}, \epsilon_{2}}$, \be \mathbb{P} \left( X \not\in
[(1-\epsilon_{1})f(N), (1+\epsilon_{1})f(N)] \right) < \epsilon_{2}.
\ee In particular we shall use this notation when $X$ is just a
function of $N$ (hence not $\lq$random'). In practice $X$ will in
this case be the expectation of some other random variable.

By $f(x) = O(g(x))$ we mean that there exist
constants $x_0$ and $C$ such that for all $x \ge x_0$, $|f(x)| \le
Cg(x)$. By $f(x) = \Theta(g(x))$ we mean that both $f(x) = O(g(x))$ and
$g(x) = O(f(x))$ hold. Finally, if
$\lim_{x\to\infty} f(x) / g(x) = 0$ then we write $f(x) = o(g(x))$.
\\
\\
Our main findings can be summed up in the following theorem.

\begin{thm}\label{thm:mainuniform}
Let $p : \mathbb{N} \rightarrow (0,1)$ be any function such that
\be\label{eq:old13} N^{-1}\ =\ o(p(N)) \ \ \ \ {\rm and}\ \ \ \
p(N)\ =\ o(1).\ee For each $N \in \mathbb{N}$ let $A$ be a random
subset of $I_{N}$ chosen according to a binomial distribution with
parameter $p(N)$. Then, as $N \rightarrow \infty$, the probability
that $A$ is difference dominated tends to one.
\par More precisely, let $\mathscr{S},
\mathscr{D}$ denote respectively the random variables $|A+A|$ and $|A-A|$.
Then the following three situations arise :
\\
\\
(i) $p(N) = o(N^{-1/2})$ : Then \be\label{eq:old14} \mathscr{S}\ \sim\
{(N\cdot p(N))^{2} \over 2} \;\;\; {\hbox{and}} \;\;\; \mathscr{D}
\sim 2\mathscr{S}\ \sim\ (N \cdot p(N))^{2}. \ee (ii) $p(N) = c \cdot
N^{-1/2}$ for some $c \in (0,\infty)$ : Define the function $g :
(0,\infty) \rightarrow (0,2)$ by \be g(x)\ :=\ 2\left(\frac{e^{-x} -
(1-x)}{x}\right). \ee Then \be\label{eq:old16} \mathscr{S}\ \sim\
g\left({c^{2} \over 2}\right) N \;\;\; {\hbox{and}} \;\;\;
\mathscr{D}\ \sim\ g(c^{2}) N. \ee (iii) $N^{-1/2} = o(p(N))$ : Let
$\mathscr{S}^{c} := (2N+1) - \mathscr{S}$, $\mathscr{D}^{c} :=
(2N+1) - \mathscr{D}$. Then \be \mathscr{S}^{c}\ \sim\ 2 \cdot
\mathscr{D}^{c}\ \sim\ {4 \over p(N)^{2}}. \ee
\end{thm}

\begin{rek}
Obviously, not all functions $p : \mathbb{N} \rightarrow (0,1)$ satisfying (1.3)
conform to the requirements of (i), (ii) or (iii) above, but these are the natural
functions to investigate in the current context. Similar remarks apply to Theorem
3.1 and Conjecture 4.2 below.
\end{rek}

Theorem \ref{thm:mainuniform} proves the conjecture in \cite{MO} and
re-establishes the validity of Nathanson's claim in a broad setting.
It also identifies the function $N^{-1/2}$ as a \emph{threshold
function}, in the sense of \cite{JLR}, for the ratio of the size of
the difference- to the sumset for a random set $A \subseteq I_{N}$.
Below the threshold, this ratio is almost surely $2 + o(1)$, above
it almost surely $1 + o(1)$. Part (ii) tells us that the ratio
decreases continuously (a.s.) as the threshold is crossed. Below the
threshold, part (i) says that most sets are $\lq$nearly Sidon sets',
that is, most pairs of elements generate distinct sums and
differences. Above the threshold, most numbers which can be in the
sumset (resp., difference set) usually are, and in fact most of these
in turn have many different representations as a sum (resp., a
difference). However the sumset is usually missing about twice as
many elements as the difference set. Thus if we replace $\lq$sums'
(resp., $\lq$differences') by $\lq$missing sums' (resp., $\lq$missing
differences'), then there is still a symmetry between what happens
on both sides of the threshold.
\par We prove Theorem \ref{thm:mainuniform}
in the next section. Our strategy will consist of first establishing
an estimate for the expectation of the random variable $\mathscr{S}$
or $\mathscr{D}$, followed by establishing sufficiently strong
concentration of these variables about their mean values. For the
second part of this strategy we will use different approaches for $p
= O(N^{-1/2})$ and $N^{-1/2} = o(p(N))$. In the former range a
fairly straightforward second moment argument works. In the latter
range, however, we will employ a specialization of
the Kim-Vu martingale lemma (Lemma 3.1 in
\cite{Vu2}, Lemma 2.2 below).
\\
\\
In Section 3, we extend our result to arbitrary binary linear forms.
The paper \cite{NOORS} provides motivation for studying these
objects. By a \emph{binary linear form} we mean a function $f(x,y) =
ux+vy$ where $u,v \in \mathbb{Z}_{\neq 0}$, $u \geq |v|$ and
GCD$(u,v) = 1$. For a set $A$ of integers we let \be f(A) := \{
ua_{1}+va_{2} : a_{1},a_{2} \in A \}. \ee Except in the special case
$u=v=1$ we always have that $f(x,y) \neq f(y,x)$ whenever $x \neq
y$. Thus we refer to $f$ as a \emph{difference form} and a set
$f(A)$ as a \emph{generalized difference set}, whenever $u > |v|$.
Theorem \ref{thm:binaryforms} allows us to compare the sizes of
$f(A)$ and $g(A)$ for random sets $A$, and arbitrary difference
forms $f$ and $g$ when $N^{-3/5} = o(p(N))$.
\\
\\
Two situations arise :
\\
\\
(a) for some pairs of forms, the same one a.s. dominates the
other for all parameters $p=p(N)$ in this range.
In fact every other difference form
dominates $x-y$, and hence also $x+y$.
\\
\\
(b) for certain pairs $f$ and $g$, something very nice happens.
Namely, there is now a \emph{sharp threshold}, in the sense of
\cite{JLR}, at $c_{f,g} N^{-1/2}$, for some computable constant
$c_{f,g} > 0$, depending on $f$ and $g$. One form dominates a.s.
below the threshold, and the other one a.s. above it. This fact may
be considered a partial generalization of the main result of
\cite{NOORS} to random sets, partial in the sense that is only
applies to certain pairs of forms. Namely, they proved that for any
two forms $f$ and $g$ (including $x+y$), there exist finite sets
$A_{1}$ and $A_{2}$ such that $|f(A_{1})| > |g(A_{1})|$ whereas
$|f(A_{2})| < |g(A_{2})|$.
\par We leave it to future work to investigate what happens as the threshold
is crossed in this situation.
\\
\\
In Section 4 we make a brief summary of this and other remaining
questions and make suggestions for other problems to study. In
particular, we suggest looking at other probabilistic models for
choosing random sets. This is partly motivated by the fact that our
results in Section 3 only apply when $N^{-3/5}=o(p(N))$. The reason
is that, for faster decaying $p(N)$, as we shall see, the variance
in the size of the random set $A$ itself swamps all other error
terms, and it is meaningless to compare $|f(A)|$ and $|g(A)|$; in other words, the
model itself becomes useless. This may be considered a problem
when $N^{-3/4}=o(p(N))$. For $p(N) = o(N^{-3/4})$, the results of
\cite{GJLR} imply that all pairs $(x,y)$ in a random set a.s.
generate different values $f(x,y)$, for any $f$, so that $|f(A)| =
|g(A)|$ a.s. for any $f$ and $g$.

\setcounter{equation}{0}
\section{Proof of Theorem \ref{thm:mainuniform}}\label{sec:generalizations12}

Our strategy for proving the various assertions in Theorem
\ref{thm:mainuniform} is the following. Let $\mathscr{X}$ be one of
the random variables $\mathscr{S}, \mathscr{D}, \mathscr{S}^{c},
\mathscr{D}^{c}$, as appropriate. We then carry out the following
two steps :
\\
\\
\emph{Step 1} : Prove that $\mathbb{E}[\mathscr{X}]$ behaves asymptotically
as asserted in the theorem.
\\
\emph{Step 2} : Prove that $\mathscr{X}$ is strongly concentrated
about its mean.
\\
\\
As already mentioned, the calculations required to perform these two
steps differ according as to whether $p(N) = O( N^{-1/2})$ or
$N^{-1/2} = o(p(N))$. In particular, in the former case, \emph{Step
2} is achieved by a fairly straightforward second moment argument,
whereas a more sophisticated concentration inequality is used in the
latter case. We thus divide the proof of the theorem into two
separate cases, depending on the parameter function $p$.
\\
\\
\emph{Throughout the paper we often abuse notation to save space,
writing $p$ for $p(N)$.} As we never consider the case where $p(N)$
is constant (as this case has been analyzed in \cite{MO}), this
should not cause any confusion.
\\
\\
\textbf{Case I :} $p(N) = O(N^{-1/2})$.
\\
\\
We first concentrate on the sumset and prove the various assertions
in parts (i) and (ii) of the theorem. The proofs for the difference
set will be similar. For any finite set $A \subseteq \mathbb{N}_{0}$
and any integer $k \geq 1$, let \be A_{k}\ :=\ \left\{ \{
\{a_{1},a_{2}\},\dots,\{a_{2k-1},a_{2k}\} \} : a_{1}+a_{2} = \dots =
a_{2k-1} + a_{2k} \right\}. \ee In words, $A_{k}$ consists of all
unordered $k$-tuples of unordered pairs of elements of $A$ having
the same sum. Let $X_{k} := |A_{k}|$. So if $A$ is a random set,
then each $X_{k}$ is a non-negative integer valued random variable.
The crucial observation for our work is that, in the model we are
considering, the random variables $X_{k}$ are all highly
concentrated :

\begin{lem}\label{lem:conc-bad-tuples} For $p(N) = O(N^{-1/2})$ we have for
every $k$ that \be\label{eq:old22} \mathbb{E}[X_{k}]\ \sim\ {2 \over
(k+1)!} \left( {p(N)^{2} \over 2} \right)^{k} N^{k+1} \ee and, more
significantly, $X_{k} \sim \mathbb{E}[X_{k}]$ whenever
$N^{-\left(\frac{k+1}{2k}\right)}=o(p(N))$.
\end{lem}

\begin{proof} We write $p$ for $p(N)$. By the central limit theorem it is clear that
\be
|A|\ \sim\ Np.
\ee
Each $X_{k}$ can be written as a sum of indicator variables $Y_{\alpha}$,
one for each unordered $k$-tuple $\alpha$. There are two types of
$k$-tuples : those consisting of $2k$ distinct elements of $I_{N}$ and those
in which one element is repeated twice in one of the $k$ pairs, and the sum
of each of the $k$ pairs is even. The probability of any $k$-tuple of
the former type occurring in $A_{k}$ is $p^{2k}$, whereas for
$k$-tuples of the latter type this probability is $p^{2k-1}$.
Let there be a total of $\xi_{1,k}(N)$ $k$-tuples of the
former type and $\xi_{2,k}(N)$ of the latter type. Then, by linearity of
expectation,
\be
\mathbb{E}[X_{k}]\ =\ \xi_{1,k}(N) \cdot p^{2k} +
\xi_{2,k}(N) \cdot p^{2k-1}.
\ee
We have
\be
\xi_{1,k}(N)\ =\ \sum_{n=0}^{2N} \left( \begin{array}{c} R(n) \\ k
\end{array} \right),
\ee where $R(n)$ is the number of representations of $n$ as a sum of
two distinct elements of $I_{N}$, and hence we easily estimate \be
\xi_{1,k}(N)\ =\ \sum_{n=2k}^{2N-2k} \left( \begin{array}{c} \min \{
\lfloor \frac{n}{2} \rfloor, \lfloor \frac{2N-n}{2} \rfloor \} \nonumber\\
k
\end{array} \right)\ \sim\ 2 \cdot \sum_{n=2k}^{N} \left( \begin{array}{c}
\lfloor \frac{n}{2} \rfloor \\ k \end{array} \right) \ee \be\ \sim\ 2
\cdot 2 \cdot \sum_{n=1}^{\lfloor N/2 \rfloor} \left(
\begin{array}{c} n \\ k
\end{array} \right)\ \sim\ 4 \left( \begin{array}{c} \lceil N/2 \rceil \\
k+1 \end{array} \right)\ \sim\ 4 {\left( N/2 \right)^{k+1} \over
(k+1)!}.  \ee Thus \be \xi_{1,k}(N) \cdot p^{2k}\ \sim\ {2 \over
(k+1)!} \left( {p^{2} \over 2} \right)^{k} N^{k+1}. \ee A similar
calculation shows that $\xi_{2,k}(N) = O_{k}(N^{k})$, hence
$\xi_{2,k}(N) \cdot p^{2k-1} = O_{k} (N^{k}p^{2k-1})$. Since $N^{-1}
= o(p)$ it follows that \be \mathbb{E}[X_{k}] \sim {2 \over (k+1)!}
\left( {p^{2} \over 2} \right)^{k} N^{k+1}, \ee in accordance with
the lemma. To complete the proof of the lemma, we need to show that,
whenever $N^{-\left(\frac{k+1}{2k}\right)} = o(p)$, the random
variable $X_{k}$ becomes highly concentrated about its mean as $N
\rightarrow \infty$. We apply a standard second moment method. In
the notation of \cite{AS}, Chapter 4, since we already know in this
case that $\mathbb{E}[X_{k}] \rightarrow \infty$, it suffices to
show that $\Delta = o(\mathbb{E}[X_{k}]^{2}) =
o_{k}(N^{2k+2}p^{4k})$, where \be \Delta\ =\ \sum_{\alpha \sim \beta}
\mathbb{P}(Y_{\alpha} \wedge Y_{\beta}), \ee the sum being over
pairs of $k$-tuples which have at least one number in common. It is
easy to see that, since $N^{-1} = o(p)$, the main contribution to
$\Delta$ comes from pairs $\{\alpha,\beta\}$ of $k$-tuples, each of
which consist of $2k$ distinct elements of $I_{N}$, and which have
exactly one element in common. The number of such pairs is
$O_{k}(N^{2k+1})$ since there are \bi \item $O_{k}(N^{k+1})$ choices
for $\alpha$,
\item $2k$ choices for the common element with $\beta$, \item $O(N)$
choices for the sum of each pair in $\beta$,
\item $O_{k}(N^{k-1})$ choices for the remaining elements in
$\beta$. \ei
Since a total of $4k-1$ elements of $I_{N}$ occur in total in
$\alpha \cup \beta$, we have $\mathbb{P}(Y_{\alpha} \wedge
Y_{\beta}) = p^{4k-1}$. Thus \be \Delta\ =\ O_{c,k}(N^{2k+1}p^{4k-1})\
=\ o_{c,k}(N^{2k+2}p^{4k}), \;\; {\hbox{since $N^{-1} = o(p)$}}.
 \ee This completes the proof of Lemma 2.1.
\end{proof}

We can now prove parts (i) and (ii) of the theorem. First suppose $p
= o(N^{-1/2})$. By \eqref{eq:old22} we have $X_{1} \sim {1 \over 2}
N^{2}p^{2}$, whereas $X_{2} \sim {1 \over 12} N^{3}p^{4}$ for
$N^{-3/4} = o(p)$ and $\mathbb{E}[X_{2}] = O(1)$ otherwise.
\par Since $p = o(N^{-1/2})$ we have $\max(1,N^{3}p^{4}) = o(N^{2}p^{2})$ and thus
$X_{2} = o(X_{1})$ almost surely. In other words, as $N \rightarrow
\infty$, all but a vanishing proportion of pairs of element of $A$
will have distinct sums. It follows immediately that \be \mathscr{S}
\ \sim\ X_{1}\ \sim\ {(Np)^{2} \over 2} , \ee as claimed.
\\
\\
Now suppose $p = cN^{-1/2}$ for some fixed $c > 0$. This time we
will need to consider all the $X_{k}$ together. Let $\mathscr{P}$ be
the partition on $A_{1}$ whereby $\{a_{1},a_{2}\}$ and
$\{a_{3},a_{4}\}$ are in the same part if and only if $a_{1}+a_{2} =
a_{3}+a_{4}$. For each $i > 0$ let $\tau_{i}$ denote the number of
parts of size $i$ (as a random variable). Observe that \be
\mathscr{S}\ =\ \sum_{i=0}^{\infty} \tau_{i} \ee and, for each $k \geq
1$, that \be\label{eq:old210} \sum_{i=1}^{\infty} \left( \begin{array}{c} i \\
k
\end{array} \right) \tau_{i}\ =\ X_{k}. \ee \eqref{eq:old210} is a system of
infinitely many equations in the variables $\tau_{i}$, which
together determine $\mathscr{S}$. For any $m \geq 1$, the general
solution of the subsystem formed by the first $m$ equations (i.e.:
$k=1,\dots,m$) is readily checked to be \be\label{eq:old211}
\mathscr{S}\ =\ \sum_{k=1}^{m} (-1)^{k-1} X_{k} +
\sum_{k=m+1}^{\infty} \left\{ \sum_{i=0}^{m} (-1)^{i} \left(
\begin{array}{c} k \\ i \end{array} \right) \right\} \tau_{k}. \ee
Regarding the second sum on the right of \eqref{eq:old211} we have
that \be \left| \sum_{k=m+1}^{\infty} \left\{ \sum_{i=0}^{m}
(-1)^{i} \left(
\begin{array}{c} k \\ i \end{array} \right) \right\} \tau_{k}
\right|\ \leq\ \sum_{k=m+1}^{\infty} \left( \begin{array}{c} k \\ m
\end{array} \right) \tau_{k}\ =\ X_{m} - \tau_{m} \leq X_{m}. \ee
Hence it follows that, for any $m \geq 1$, \be \left| X -
\sum_{k=1}^{m} (-1)^{k-1} X_{k} \right|\ \leq\ X_{m}. \ee Now for $p =
cN^{-1/2}$, Lemma 2.1 says that \be X_{m}\ \sim\ 2 {\left( {c^{2}
\over 2} \right)^{m} \over (m+1)!} N,  \ee and since ${\left( {c^{2}
/ 2} \right)^{m} \over (m+1)!} \rightarrow 0$ as $m \rightarrow
\infty$, another application of Lemma 2.1 implies that \be
\mathscr{S}\ \sim\ \sum_{k=1}^{\infty} (-1)^{k-1} X_{k}\ \sim\ 2 \cdot
\left( \sum_{k=1}^{\infty} {(-1)^{k-1} \left( \frac{c^{2}}{2}
\right)^{k} \over (k+1)!} \right) \cdot N. \ee So to prove
\eqref{eq:old16} it just remains to verify that \be g(x)\ =\ 2 \cdot
\sum_{k=1}^{\infty} {(-1)^{k-1} x^{k} \over (k+1)!},  \ee which is
an easy exercise.

This proves parts (i) and (ii) of Theorem 1.1 for the sumset. For
the difference set one reasons in an entirely parallel manner. One
now defines, for each $k \geq 1$, \be A^{\prime}_{k} := \left\{ \{
(a_{1},a_{2}),\dots,(a_{2k-1},a_{2k}) \} : a_{1} - a_{2} = \dots =
a_{2k-1} - a_{2k} \neq 0 \right\}. \ee In words, $A^{\prime}_{k}$
consists of all $k$-tuples of ordered pairs of elements of $A$ which
have the same non-zero difference. We let $X^{\prime}_{k} :=
|A_{k}|$ and in a completely analogous manner to Lemma 2.1 prove
that \be\label{eq:old216} \mathbb{E}[X^{\prime}_{k}]\ \sim\ {2 \over
(k+1)!} p^{2k} N^{k+1}, \ee and that $X^{\prime}_{k} \sim
\mathbb{E}[X^{\prime}_{k}]$ whenever
$N^{-\left(\frac{k+1}{2k}\right)} = o(p)$. We define the partition
$\mathscr{P}^{\prime}$ of $A^{\prime}_{1}$ in the obvious way and
let $\tau^{\prime}_{i}$ denote the number of parts of size $i$, for
each $i \geq 1$. Since $\mathscr{D} = 1 + \sum_{i = 1}^{\infty}
\tau_{i}$ we can follow exactly the same analysis as above to deduce
\eqref{eq:old14} and \eqref{eq:old16}. This completes the proofs of
parts (i) and (ii) of the theorem.
\\
\\
\textbf{Case II :} $N^{-1/2} = o(p(N))$.
\\
\\
Recall $p(N) = o(1)$. Set $p = p(N)$ and $P = 1/p$; thus $P =
o(N^{1/2})$ (as $p = o(1)$ we have $\lim_{N\to\infty} P = \infty$).
Again we begin with the sumset. Recall the two steps to be
accomplished :
\\
\\
\emph{Step 1} : We prove that $\mathbb{E}[\mathscr{S}^{c}] \sim 4P^{2}$.
\\
\emph{Step 2} : We prove that the random variable $\mathscr{S}^{c}$ is
strongly concentrated about its mean.
\\
\\
We begin with the simpler \emph{Step 1}. For each $n \in I_{2N}$,
let $\mathscr{E}_{n}$ denote the event that $n \not\in A+A$. Thus
\be \mathbb{E}[\mathscr{S}^{c}] = \sum_{n=0}^{2N}
\mathbb{P}(\mathscr{E}_{n}).  \ee Observe that
$\mathbb{P}(\mathscr{E}_{n}) = \mathbb{P}(\mathscr{E}_{2N-n})$.
Since all the ways of representing any given $n$ as a sum of two
elements of $I_{N}$ are independent of one another, we have, for $n
\in I_{N}$, \be\label{eq:old217} \twocase{
\mathbb{P}(\mathscr{E}_{n}) \ = \ }{(1-p^{2})^{n/2} (1-p)}{if $n$ is
even}{(1-p^{2})^{(n+1)/2}}{if $n$ is odd.} \ee Since $p = o(1)$ we
have $1-p \sim 1$, and since $N^{-1/2} = o(p)$ we have
$(1-p^{2})^{N} = o(1)$. Thus it is easy to see that
\be\label{eq:expvalueESc} \mathbb{E}[\mathscr{S}^{c}]\ \sim\ 4 \cdot
\sum_{m=0}^{\lfloor N/2 \rfloor} (1-p^{2})^{m}\ \sim\ {4 \over p^{2}}\
=\ 4P^{2},  \ee as claimed. This completes \emph{Step 1}.
\\
\\
For \emph{Step 2} we need the martingale machinery of Kim and Vu.
\\
\\
We use notation consistent with \cite{Vu2}. Consider a fixed $N$,
which shall tend to infinity in our estimates. Let $\Omega :=
\{0,1\}^{N+1}$. Thus every subset $A$ of $I_{N}$ can be identified
with an element of $\Omega$. We are working in the probability space
$(\Omega,\mu)$ where $\mu$ is the product measure with parameter
$p$. For each $A \in \Omega$, $n \in I_{N}$ and $x \in \{0,1\}$,
define \be\label{eq:defnCnxA} C_{n}(x,A)\ :=\ \left|
\mathbb{E}\left[\mathscr{S}^{c}|a_{0},...,a_{n-1},a_{n}=x\right] -
\mathbb{E}\left[\mathscr{S}^{c}|a_{0},...,a_{n-1}\right] \right|;
\ee by $\mathbb{E}\left[\mathscr{S}^{c}|a_{0},...,a_{m}\right]$ we
mean the expected value of the random variable $\mathscr{S}^{c}$,
given that for $k \in \{0,\dots,m\}$ the element $k$ is always
(resp., never) in the subset if $a_k = 1$ (resp., $a_k=-1$). Let \be
C(A)\ :=\ \max_{n,x} C_{n}(x,A). \ee Further put \be\label{eq:old220}
V_{n}(A)\ :=\ \int_{0}^{1} C_{n}^{2}(x,A) d^{n}\mu = pC_{n}^{2}(0,A) +
(1-p)C_{n}^{2}(1,A) \ee and \be\label{eq:old221} V(A)\ :=\
\sum_{n=0}^{N} V_{n}(A). \ee For two arbitrary positive numbers
$\mathbf{V}$ and $\mathbf{C}$, define the event \be
\mathbb{B}_{\mathbf{V,C}}\ :=\ \{A : C(A) \geq \mathbf{C} \;
{\hbox{or}} \; V(A)\ \geq\ \mathbf{V} \}. \ee Then the following is a
specialization of a result appearing in \cite{Vu2} :

\begin{lem}\label{lem:martingaleconc}
For any positive numbers $\lambda,\mathbf{V},\mathbf{C}$ such that
$\lambda \leq 4\mathbf{V}/\mathbf{C}^{2}$ we have
\be\label{eq:old223} \mathbb{P}\left( \left| \mathscr{S}^{c} -
\mathbb{E}[\mathscr{S}^{c}] \right| \geq \sqrt{\lambda \mathbf{V}}
\right)\ \leq\ 2e^{-\lambda/4} +
\mathbb{P}(\mathbb{B}_{\mathbf{V,C}}). \ee
\end{lem}

We quickly sketch how Lemma \ref{lem:martingaleconc} completes the
proof of assertion (iii) of Theorem \ref{thm:mainuniform}.  We shall
take \be  \;\;\; \lambda := \kappa_0\log P, \;\;\; \mathbf{V} :=
\kappa_{1}(P \log P)^{3}, \;\;\; \mathbf{C} := \kappa_{2}P \log P.
\ee We show that for appropriately chosen $\kappa_{1},\kappa_{2}$ we
have \be\label{eq:probvco1} \mathbb{P}(\mathbb{B}_{\mathbf{V,C}})\ =\
o(1). \ee From \eqref{eq:old223} and $\lim_{N\to\infty} P = \infty$,
for sufficiently small $\kappa_0$ we will then be able to conclude
that \be |\mathscr{S}^{c} - \mathbb{E}[\mathscr{S}^{c}]|\ =\ O(P^{3/2}
\log^2 P) \;\; {\hbox{a.s. as $N \rightarrow \infty$}}. \ee As
$\mathbb{E}[\mathscr{S}^{c}] \sim 4P^2$ (see
\eqref{eq:expvalueESc}), assertion (iii) in Theorem
\ref{thm:mainuniform} follows immediately. Thus we are reduced to
proving \eqref{eq:probvco1}, which we now proceed to do.
\\
\\
First we simplify things a little. For any $n \in I_{N}$ and $A \in
\Omega$, we introduce the shorthand \be \mathscr{U}_{n,A}\ :=\
\mathscr{S}^{c}|a_{0},...,a_{n-1}. \ee Let \be \Delta_{n}(A)\ :=\
\mathbb{E}[\mathscr{U}_{n,A}|a_{n}=0] -
\mathbb{E}[\mathscr{U}_{n,A}|a_{n}=1]. \ee As
$\mathbb{E}[\mathscr{U}_{n,A}|a_{n}=0] \geq
\mathbb{E}[\mathscr{U}_{n,A}|a_{n}=1]$, we see $\Delta_{n}(A) \ge
0$. For $x \in \{0,1\}$, \be C_{n}(x,A)\ =\
|\mathbb{E}[\mathscr{U}_{n,A}|a_{n}=x] -
\mathbb{E}[\mathscr{U}_{n,A}]|. \ee Since \bea
\mathbb{E}[\mathscr{U}_{n,A}]&\ =\ &
p\mathbb{E}[\mathscr{U}_{n,A}|a_{n}=1] +
(1-p)\mathbb{E}[\mathscr{U}_{n,A}|a_{n}=0] \nonumber\\ & \ = \ &
\mathbb{E}[\mathscr{U}_{n,A}|a_{n}=1] + (1-p) \Delta_n(A),
 \eea $C_n(x,A)$ from \eqref{eq:defnCnxA} simplifies to \be\label{eq:old227}
 C_{n}(1,A)\ =\ (1-p)\Delta_{n}(A), \;\;\;
C_{n}(0,A)\ =\ p\Delta_{n}(A). \ee Since $p < 1-p$ for all
sufficiently large $N$, we have then \be C(A) = (1-p) \max_{0 \leq n
\leq N} \Delta_{n}(A). \ee Further, \eqref{eq:old220},
\eqref{eq:old221} and \eqref{eq:old227} yield \be\label{eq:old230}
V(A) = p(1-p) \sum_{n=0}^{N} \Delta_{n}^{2}(A). \ee This completes
our simplifications.

Recall that in order to use \eqref{eq:old223} from Lemma
\ref{lem:martingaleconc} we need to prove \eqref{eq:probvco1}
(namely that $\mathbb{P}(\mathbb{B}_{\mathbf{V,C}}) = o(1)$). The
heart of the proof of \eqref{eq:probvco1} is to show that for an
appropriate choice of $\kappa_{1}, \kappa_{2}, \kappa_{3} > 0$, with
probability $1- o(1)$ all three of the following events occur:
\be\label{eq:old231} \sum_{|N-n|
> \kappa_{3} P^{2} \log P} \Delta_{n}(A) = o(1), \ee \be\label{eq:old232}
\max_{0 \leq |N-n| \leq \kappa_{3} P^{2} \log P} \Delta_{n}(A) \leq
\kappa_{2} P \log P, \ee \be\label{eq:old233} V(A) \leq \kappa_{1}
(P \log P)^{3}. \ee We claim that
\eqref{eq:old231}-\eqref{eq:old233} imply \eqref{eq:probvco1}. This
follows immediately from applying the trivial bound
$\mathbb{P}(\mathbb{B}_{\mathbf{V,C}}) \le \mathbb{P}(V \ge
\mathbf{V}) + \mathbb{P}(C \ge \mathbf{C})$ and using
\eqref{eq:old231}-\eqref{eq:old233} to show these two probabilities
are both $o(1)$.

To summarize, the proof is completed by verifying
\eqref{eq:old231}-\eqref{eq:old233}. Notice also that
\eqref{eq:old231} and \eqref{eq:old232}, together with
\eqref{eq:old230}, imply \eqref{eq:old233} for any choice of
$\kappa_{1} > \kappa_{2}^2 \kappa_{3}$, so it just remains to prove
the former two. As in the arguments that follow there is a symmetry
between $n$ and $N-n$, we consider $n$ with $0 \leq n \leq N/2$; the
remaining $n$ follow similarly.
\\
\\
First, consider \eqref{eq:old231}. Note that, depending on the
parameter $P$, this sum could be empty. This will not affect the
argument to follow. The proof is by an averaging argument, i.e.: for
each $n$ we first consider
$\mathbb{E}_{a_{0},...,a_{n-1}}[\Delta_{n}(A)]$. This quantity has a
very natural interpretation : in words, it is the expected increase
in the size of the sumset $A+A$ brought about by the addition of the
number $n$ to $A$. For every $z \in \{n,...,n+N\}$, adding $n$ to
$A$ will add $z$ to $A+A$ if and only if $z-n \in A$ and, for any
other numbers $n_{1},n_{2}$ such that $n_{1}+n_{2} = z$, either
$n_{1} \not\in A$ or $n_{2} \not\in A$. Let $\mathscr{E}_{z}$ be the
event that $z$ gets added to $A+A$ by the addition of $n$ to $A$.
Then using \eqref{eq:old217} we can explicitly estimate \be
\mathbb{E}_{a_{0},...,a_{n-1}}[\Delta_{n}(A)]\ =\ \sum_{z=n}^{n+N}
\mathbb{P}[\mathscr{E}_{z}]\ \sim\ p \sum_{z=n}^{n+N} (1-p^{2})^{\min
\{ \lfloor \frac{z}{2} \rfloor, \lfloor \frac{2N-z}{2} \rfloor \}}.
\ee Since $n \leq N/2$, the last sum is asympotically no more than
\be 2p\sum_{r=\lfloor n/2 \rfloor}^{\lfloor N/2 \rfloor}
(1-p^{2})^{r}\ =\ (2+o(1)) P (1-p^{2})^{n/2}.  \ee By Markov's
inequality, we deduce that for any $n$, \be \mathbb{P}(\Delta_{n}(A)
\geq 2P(1-p^{2})^{n/4})\ \leq\ (1+o(1))(1-p^{2})^{n/4}. \ee Then, just
using a trivial union bound \be \mathbb{P} \left( \bigvee
\mathscr{E}_{n} \right)\ \leq\ \sum \mathbb{P}(\mathscr{E}_{n}),  \ee
it follows that, with probability at least \be 1 -
\left(1+o(1)\right) \cdot \sum_{n = \lceil \kappa_{3} P^{2} \log P
\rceil}^{\lfloor N/2 \rfloor} (1-p^{2})^{n/4},  \ee we have \be
\Delta_{n}(A)\ \leq\ 2P(1-p^{2})^{n/4} \;\; {\hbox{for all $n$ such
that $\kappa_{3} P^{2} \log P \leq n \leq N/2$}}.  \ee Then
\eqref{eq:old231} clearly follows provided \be P \cdot \sum_{n =
\lceil \kappa_{3} P^{2} \log P \rceil}^{\lfloor N/2 \rfloor}
(1-p^{2})^{n/4}\ =\ o(1), \ee which is clearly the case for
sufficiently large $\kappa_{3}$, since $1=o(P)$.
\\
\\
We now turn to \eqref{eq:old232}. Firstly, a similar argument to the
one just given shows that, even if $n \leq \kappa_{3} P^{2} \log P$,
adding $n$ to a random set $A$ is very probably not going to add any
elements at all to $A+A$ which are larger than $\kappa_{3} P^{2}
\log P$. Secondly, among the numbers in $I_{\kappa_{3} P^{2} \log
P}$, the addition to $A$ of one number cannot add to $A+A$ more
numbers than were in $A$ already, plus maybe one more. But
Chernoff's inequality (\cite{AS}, Corollary A.14) implies that, with
probability $1-e^{-c_{1}\kappa_{3}P\log P}$, where $c_{1}$ is some
universal positive constant, $|A \cap I_{\kappa_{3} P^{2} \log P}|
\leq 2\kappa_{3} P \log P$. Then \eqref{eq:old232} follows from a
simple union bound, as long as $\kappa_{2} > 2\kappa_{3}$ for
example.
\par This completes the proof of the assertion of Theorem
\ref{thm:mainuniform}(iii)
as regards the sumset.
\\
\\
For the difference set, we proceed in two identical steps. First
consider the estimate of $\mathbb{E}[\mathscr{D}^{c}]$. Let
$\mathscr{E}_{n}$ now denote instead the event that $n \not\in A-A$
for each $n \in \pm I_{N}$. Clearly, \be \mathbb{E}[\mathscr{D}^{c}]
\ =\ 2 \cdot \sum_{n=1}^{N} \mathbb{P}(\mathscr{E}_{n}) + o(1).  \ee
For each $n > 0$ we have \be \mathscr{E}_{n}\ =\ \bigwedge_{m=0}^{N-n}
\overline{\mathscr{B}}_{m,n},  \ee where $\mathscr{B}_{m,n}$ is the
(bad) event that both $m$ and $m+n$ lie in $A$ and
$\overline{\mathscr{B}}_{m,n}$ is the complementary event. These
events are not independent, but the dependencies will not affect our
estimates. To see this rigorously, one can for example use Janson's
inequality (see \cite{AS}, Chapter 8, though this is certainly
overkill!) \be M\ \leq\ \mathbb{P} \left( \bigwedge_{m}
\overline{\mathscr{B}}_{m,n} \right)\ \leq\ M \exp \left( {\Delta
\over 1-\epsilon} \right),
 \ee where all $\mathbb{P}(\mathscr{B}_{m,n}) \leq
\epsilon$, $M = \prod_{m} \mathbb{P}(\overline{\mathscr{B}}_{m,n})$
and \be \Delta = \sum_{m \sim m^{\prime}}
\mathbb{P}(\mathscr{B}_{m,n} \wedge \mathscr{B}_{m^{\prime},n}),
 \ee the sum being over dependent pairs $\{m,m^{\prime}\}$,
i.e.: pairs such that $m^{\prime} = m+n$.
\par Note that we can take
$\epsilon = p^{2}$, we have $M = (1-p^{2})^{N-n+1}$ and \be
\twocase{\Delta\ = \ }{(N-2n+1)p^{3}}{if $n \le N/2$}{0}{if $n > N/2$,}
\ee since there is a 1-1 correspondence between dependent pairs and
3-term arithmetic progressions in $I_{N}$ of common difference $n$.
It is then easy to see that this correction term can be ignored when
we make the estimate \be \mathbb{E}[\mathscr{D}^{c}]\ \sim\ 2 \cdot
\sum_{n=1}^{N} (1-p^{2})^{N-n+1}\ \sim\ {2 \over p^{2}}\ \sim\ {1 \over
2} \mathbb{E}[\mathscr{S}^{c}],  \ee as desired.
\par The concentration of $\mathscr{D}^{c}$ about its mean can be established
in the same way as we did with $\mathscr{S}^{c}$ above. A little
more care is required in estimating quantities analogous to
$\mathbb{E}_{a_{0},...,a_{n-1}}[\Delta_{n}(A)]$, because of the
dependencies between different representations of the same
difference, but Janson's inequality can again be used to see
rigorously that this will not affect our estimates. We omit further
details and simply note that we will again obtain the result that
\be |\mathscr{D}^{c} - \mathbb{E}[\mathscr{D}^{c}]|\ =\ O(P^{3/2}
\log^{2} P) \;\; {\hbox{a.s. as $N \rightarrow \infty$}}. \ee This
completes the proof of Theorem \ref{thm:mainuniform}. \hfill $\Box$


\setcounter{equation}{0}
\section{General Binary Linear Forms}\label{sec:binaryforms}

We have the following generalization of Theorem
\ref{thm:mainuniform} :

\begin{thm}\label{thm:binaryforms}

Let $p:\N \to (0,1)$ be a function satisfying \eqref{eq:old13}. Let
$u,v$ be non-zero integers with $u \geq |v|$, GCD$(u,v) = 1$ and
$(u,v) \neq (1,1)$. Put $f(x,y) := ux+vy$. For a positive integer
$N$, let $A$ be a random subset of $I_{N}$ obtained by choosing each
$n \in I_{N}$ independently with probability $p(N)$. Let
$\mathscr{D}_{f}$ denote the random variable $|f(A)|$. Then the
following three situations arise :
\\
\\
(i) $p(N) = o(N^{-1/2})$ : Then \be \mathscr{D}_{f} \sim (N \cdot
p(N))^{2}. \ee (ii) $p(N) = c\cdot N^{-1/2}$ for some $c \in
(0,\infty)$ : Define the function $g_{u,v} : (0,\infty) \rightarrow
(0,u+|v|)$ by \be\label{eq:old32} g_{u,v}(x)\ :=\ (u+|v|) - 2|v|
\left( {1-e^{-x} \over x} \right) - (u-|v|)e^{-x}. \ee Then
\be\label{eq:old33} \mathscr{D}_{f}\ \sim\ g_{u,v} \left( {c^{2} \over
u} \right) N. \ee (iii) $N^{-1/2} = o(p(N))$ : Let
$\mathscr{D}_{f}^{c} := (u+|v|)N - \mathscr{D}_{f}$. Then \be
\mathscr{D}_{f}^{c}\ \sim\ {2u|v| \over p(N)^{2}}. \ee
\end{thm}

\begin{proof} One follows exactly the method of proof of Theorem
\ref{thm:mainuniform}, so we only give a sketch here.
\\
\\
\textbf{Case I :} $p(N) = O(N^{-1/2})$.
\\
\\
We again write $p$ for $p(N)$. For any finite set $A \subseteq
\mathbb{N}_{0}$ and any integer $k \geq 1$, let \be A^{\prime}_{k,f}\
:=\ \left\{ \{ (a_{1},a_{2}),\dots,(a_{2k-1},a_{2k}) \} :
f(a_{1},a_{2}) = \cdots = f(a_{2k-1},a_{2k}) \right\}. \ee Let
$X^{\prime}_{k,f} := |A_{k,f}|$. Then \eqref{eq:old216} has the
following generalization : \be\label{eq:old36}
\mathbb{E}[X^{\prime}_{k,f}]\ \sim\ \left( {2|v| \over (k+1)!} +
{u-|v| \over k!} \right) p^{2k} N^{k+1}, \ee and $X^{\prime}_{k,f}
\ \sim\ \mathbb{E}[X^{\prime}_{k,f}]$ whenever
$N^{-\left(\frac{k+1}{2k}\right)}=o(p)$.
\par We
shall just sketch the proof that $\mathbb{E}[X^{\prime}_{k,f}]$
behaves like the right-hand side of \eqref{eq:old36} in the case
when $v > 0$. The proof for $v < 0$ is similar, and the
concentration of $X_{k}$ about its mean when
$N^{-\left(\frac{k+1}{2k}\right)} = o(p)$ is established by the same
kind of second moment argument as in Section 2.
\par If $v > 0$ then for any $A \subseteq I_{N}$ we have
$f(A) \subseteq I_{(u+v)N}$. Then \be\label{eq:old37}
\mathbb{E}[X^{\prime}_{k,f}]\ \sim\ \xi_{k,f}(N) \cdot p^{2k}, \ee
where \be \xi_{k,f}(N)\ =\ \sum_{n=0}^{(u+v)N} \left( \begin{array}{c}
R(n) \\ k
\end{array} \right)
\ee and $R(n)$ denotes the number of solutions to the equation
$ux+vy=n$ satisfying $(x,y) \in I_{N} \times I_{N}$. For any integer
$n$, the general integer solution to $ux+vy=n$ is of course \be
x=nx_{0}-vt, \;\;\; y=ny_{0}+ut, \;\;\; t \in \mathbb{Z}, \ee where
$ux_{0}+vy_{0}=1$. If $n > 0$ then there are $\lfloor n/uv \rfloor +
O(1)$ solutions in non-negative integers, and for all such
solutions, $(x,y) \in I_{\lfloor n/u \rfloor} \times I_{\lfloor n/v
\rfloor}$. For $n \in I_{(u+v)N}$ the following three situations
then arise :
\\
\\
(I) $n \in I_{vN}$ : then all non-negative solutions satisfy $(x,y)
\in I_{N} \times I_{N}$, so $R(n) = \lfloor n/uv \rfloor + O(1)$ in
this case.
\\
(II) $vN < n \leq uN$ : we have $R(n) = \lfloor N/u \rfloor + O(1)$
for any such $n$.
\\
(III) $uN < n \leq (u+v)N$ : we have $R(n) = \lfloor \frac{1}{uv}
\left[ (u+v)N - n\right] \rfloor + O(1)$ for these $n$.
\\
\\
Thus it follows that
\be
\xi_{k,f} \sim \sum_{n=0}^{vN} \left( \begin{array}{c} \lfloor n/uv \rfloor \\
k \end{array} \right) + \sum_{n=vN}^{uN} \left( \begin{array}{c}
\lfloor N/v \rfloor \\ k \end{array} \right) + \sum_{n=uN}^{(u+v)N}
\left( \begin{array}{c} \lfloor \frac{1}{uv} \left[ (u+v)N -
n\right] \rfloor \\ k \end{array} \right) \ee \be \sim 2 \cdot uv
\sum_{n=0}^{\lfloor N/u \rfloor} \left( \begin{array}{c} n \\ k
\end{array} \right) + (u-v)N \left( \begin{array}{c} \lfloor N/u
\rfloor \\ k \end{array} \right)  \ee \be \sim 2uv {\left( {N \over
u} \right)^{k+1} \over (k+1)!} + (u-v)N {\left( {N \over u}
\right)^{k} \over k!}  \ee which, together with \eqref{eq:old37},
verifies our claim that $\mathbb{E}[X^{\prime}_{k,f}]$ behaves like
the right-hand side of \eqref{eq:old36}.
\par Once we have \eqref{eq:old36} then, in a similar manner to Section 2, we can
prove part (i) of Theorem \ref{thm:binaryforms} by noting that
$X^{\prime}_{2,f} = o(X^{\prime}_{1,f})$ almost surely when $p =
o(N^{-1/2})$, and part (ii) by showing that \be\label{eq:old39}
\mathscr{D}_{f}\ \sim\ \sum_{k=1}^{\infty} (-1)^{k-1} X^{\prime}_{k,f}
\ee when $p = cN^{-1/2}$. It's a simple exercise to check that
\eqref{eq:old39} and \eqref{eq:old36} yield \eqref{eq:old33}.
\\
\\
\textbf{Case II :} $N^{-1/2} = o(p(N))$.
\\
\\
We give a sketch of the estimate for
$\mathbb{E}[\mathscr{D}_{f}^{c}]$, the details of the concentration
estimate being completely analogous to what has gone before. Let us
continue to assume $v > 0$, the proof for $v < 0$ being similar. As
in the proof of Theorem \ref{thm:mainuniform}(iii) one may check
that various dependencies do not affect our estimates which, using
observations (I),(II),(III) above, lead to \be
\mathbb{E}[\mathscr{D}_{f}^{c}]\sim \sum_{n=0}^{vN}
(1-p^{2})^{\lfloor n/uv \rfloor} + \sum_{n=vN}^{uN} (1-p^{2})^{N/u}
+ \sum_{n=uN}^{(u+v)N} (1-p^{2})^{\lfloor
\frac{1}{uv}[(u+v)N-n]\rfloor}  \ee \be \sim 2\cdot uv
\sum_{n=0}^{\lfloor N/u \rfloor} (1-p^{2})^{n} +
(u-v)N(1-p^{2})^{N/u}.  \ee The sum is $\sim 1/p^{2}$ and the second
term is negligible since $1=o(Np^{2})$, so
$\mathbb{E}[\mathscr{D}_{f}^{c}] \sim 2uv/p^{2}$, as claimed.
\par This completes the proof of Theorem \ref{thm:binaryforms}.
\end{proof}

As mentioned earlier,
the main result of \cite{NOORS} was that, for any two binary forms $f$ and $g$,
including the case when $g(x,y) = x+y$,
there exist finite sets $A_{1},A_{2}$ of integers such that
$|f(A_{1})| > |g(A_{1})|$
and $|f(A_{2})| < |g(A_{2})|$. Theorem \ref{thm:binaryforms} has a number of
consequences on the matter of comparing $|f(A)|$ and $|g(A)|$ for given $f$
and $g$ and random subsets $A$ of $I_{N}$ for large $N$. We now reserve
the notations
$f(x,y) := u_{1}x + v_{1}y$ and $g(x,y) := u_{2} x + v_{2}y$ for two
forms being compared. Unless
otherwise stated, we assume neither $f$ nor $g$ is the form $x+y$.
A generic form $ux+vy$ will be denoted $h(x,y)$.
\par It is convenient to formalize a piece of terminology which we used
informally in the introduction :

\begin{defi}
Let $f,g$ be two binary linear forms as above. Let $p: \N \to (0,1)$
satisfy \eqref{eq:old13}. Then we say that $f$ {\bf dominates} $g$
for the parameter $p$ if, as $N \rightarrow \infty$, $|f(A)|
> |g(A)|$ almost surely when $A$ is a random subset of $I_{N}$ obtained by

choosing each $n \in I_{N}$ independently with probability
$p(N)$.
\end{defi}

We now consider three different regimes (depending on how rapidly
$p(N)$ decays). In the arguments below we shall write $p$ for
$p(N)$. The most interesting behavior will be isolated afterwards as
Theorem \ref{thm:compareforms}.

$\;$ \\
{\bf Regime 1 :} $N^{-1/2} = o(p(N))$.
\\
\\
Then part (iii) of Theorem \ref{thm:binaryforms} implies, in
particular, that $\mathscr{D}_{f} \sim (u_{1}+|v_{1}|)N$ and
$\mathscr{D}_{g} \sim (u_{2}+|v_{2}|)N$. Hence $f$ dominates $g$
when $u_{1}+|v_{1}| > u_{2}+|v_{2}|$. In particular, this is the
case if $g(x,y) = x - y$. If $u_{1}+|v_{1}| = u_{2}+|v_{2}|$ then
the theorem says that $f$ dominates $g$ if and only if $u_{1}|v_{1}|
< u_{2}|v_{2}|$, which is the case if and only if $u_{1} > u_{2}$.
\\
\\
{\bf Regime 2 :} $p = o(N^{-1/2})$.
\\
\\
Part (i) of Theorem \ref{thm:binaryforms} says that $\mathscr{D}_{f}
\sim \mathscr{D}_{g} \sim (Np)^{2}$ for any $f$ and $g$. For every
$k \geq 1$ we have \be X^{\prime}_{k,h}\ =\ \Theta_{k,u,v}
(N^{k+1}p^{2k}). \ee Thus \be\label{eq:old310} X^{\prime}_{k+1,h}\ =\
O(Np^{2} X^{\prime}_{k,h}) = o(X^{\prime}_{k,h}) \;\; {\hbox{almost
surely.}} \ee The second moment method gives standard deviations \be
\sigma (X^{\prime}_{k,h}) = \Theta_{k,u,v} \left(
\sqrt{\mathbb{E}[X^{\prime}_{k,h}] + \Delta} \right)\ =\
\Theta_{k,u,v} \left(\max \{ N^{\frac{k+1}{2}}p^{k},
N^{k+\frac{1}{2}}p^{2k-\frac{1}{2}} \} \right). \ee In particular we
have $\sigma(X^{\prime}_{1,h}) = \Omega\left([Np]^{3/2}\right)$ and
$X^{\prime}_{2,h} = \Theta(N^{3}p^{4})$. First of all, then, if
$(N^{3}p^{4}) = O((Np)^{3/2})$, i.e.: if $p = O(N^{-3/5})$, then
the uncertainty in the size of the random set $A$ itself swamps
everything else, and our model is worthless.
\par If $N^{-3/5} = o(p)$ then, by \eqref{eq:old310}, it is in the first instance the
$X^{\prime}_{2,h}$-term which will be decisive. By \eqref{eq:old36}
we have \be X^{\prime}_{2,h}\ \sim\ {1 \over u^{2}} \left( {|v| \over
3} + {u-|v| \over 2} \right) N^{3} p^{4}. \ee Hence $f$ dominates
$g$ in this range of $\delta$ if $\alpha(u_{1},v_{1}) <
\alpha(u_{2},v_{2})$ where \be\label{eq:old313} \alpha(u,v)\ :=\ {1
\over u^{2}} \left( {|v| \over 3} + {u-|v| \over 2} \right)\ =\
{3u-|v| \over 6u^{2}}. \ee Since it is easy to see that
$\alpha(u_{1}v_{1}) = \alpha(u_{2},v_{2})$ if and only if
$u_{1}=u_{2}, v_{1} = \pm v_{2}$, this allows us to compare any pair
of forms in the range $N^{-3/5} = o(p)$ and $p = o(N^{-1/2})$,
except a pair $ux \pm vy$. But for such a pair, our methods are
entirely worthless anyway, since all the estimates in this section
depend only on $|v|$. Note in particular that $\alpha(u,v) <
\alpha(1,-1)$ for any $(u,v) \neq (1,-1)$ so that any other form
dominates $x - y$.
\\
\\
{\bf Regime 3 :} $p = cN^{-1/2}$.
\\
\\
By part (iii) of Theorem \ref{thm:binaryforms}, for a given value of
the parameter $c \in (0,\infty)$, $f$ dominates $g$ if \be
g_{u_{1},v_{1}}\left( {c^{2} \over u_{1}} \right)\ >\ g_{u_{2},v_{2}}
\left( {c^{2} \over u_{2}} \right). \ee Since $g_{u,v}(x)
\rightarrow u+|v|$ as $x \rightarrow \infty$, $f$ will dominate $g$
for sufficiently large values of $c$, provided $u_{1}+|v_{1}| >
u_{2}+|v_{2}|$. This is as expected from Regime 1. On the other hand,
the Taylor expansion of $g_{u,v}$, as a function of $c$, around $c =
0$, reads \be g_{u,v}(c)\ =\ c^{2} - \alpha(u,v) c^{4} +
O_{u,v}(c^{6}). \ee Thus $f$ dominates $g$ for sufficiently small
values of $c$ provided $\alpha(u_{1},v_{1}) < \alpha(u_{2},v_{2})$.
Again this is as expected, this time from Regime 2. Note that the
injectivity of $\alpha$ allows us to even compare forms with the
same value of $u+|v|$, namely : for a fixed value of $u+|v|$,
$\alpha(u,v)$ is clearly a decreasing function of $u$. Hence if
$u_{1}+|v_{1}| = u_{2}+|v_{2}|$ then $f$ dominates $g$ for all
values of $c \in (0,\infty)$ if and only if $u_{1} > u_{2}$. Note that this
is the same condition as in Regime 1. More
generally, we have that $f$ dominates $g$ for all values of $c$
whenever $\alpha(u_{1},v_{1}) < \alpha(u_{2},v_{2})$ and
$u_{1}+|v_{1}| \leq u_{2}+|v_{2}|$. In particular this is the case
for $g(x,y) = x-y$ and any other $f$.
\\
\\
The most interesting phenomenon arises when we compare two forms
such that \be\label{eq:old316} u_{1}+|v_{1}|\ >\ u_{2}+|v_{2}| \;\;\;
{\hbox{and}} \;\;\; \alpha(u_{1},v_{1})\ >\ \alpha(u_{2},v_{2}). \ee
Then the combined observations of Regimes 1, 2 and 3 imply that there

exists some \\ $c_{f,g} > 0$ such that \bea\label{eq:old317} & &
{\hbox{$g$ dominates $f$ whenever $N^{-3/5} = o(p)$ and $p =
o(N^{-1/2})$}} \nonumber\\ & & {\hbox{or $p = cN^{-1/2}$ for any $0
< c < c_{f,g}$}}, \eea whereas \be\label{eq:old318} {\hbox{$f$
dominates $g$ whenever $p = cN^{-1/2}$ for any $c > c_{f,g}$ or
$N^{-1/2} = o(p)$}}. \ee This observation may be considered a
partial generalization of the main result of \cite{NOORS} to random
sets, partial in the sense that it only applies to pairs of forms
satisfying \eqref{eq:old316}. Equations \eqref{eq:old317} and
\eqref{eq:old318} say that we have a sharp threshold, below which
$g$ dominates $f$ and above which $f$ dominates $g$. We leave it to
future work to determine what happens as one crosses this sharp
threshold.
\\
\\
We close this section by summarizing the most important observations
above in a theorem.

\begin{thm}\label{thm:compareforms}
Let $f(x,y) = u_{1}x+u_{2}y$ and $g(x,y) = u_{2}x+g_{2}y$, where
$u_{i} \geq |v_{i}| > 0$, GCD$(u_{i},v_{i}) = 1$ and $(u_{2},v_{2})
\neq (u_{1},\pm v_{1})$. Let $\alpha : \mathbb{Z}_{\neq 0}^{2}
\rightarrow \mathbb{Q}$ be the function given by \eqref{eq:old313}.
The following two situations can be distinguished :
\\
\\
(i) $u_{1}+|v_{1}| \geq u_{2}+|v_{2}|$ and $\alpha(u_{1},v_{1}) <
\alpha(u_{2},v_{2})$.
\\
\\
Then $f$ dominates $g$ for all $p$ such that $N^{-3/5} = o(p)$ and $
p = o(1)$. In particular, every other difference form dominates the
form $x-y$ in this range.
\\
\\
(ii) $u_{1}+|v_{1}| > u_{2}+|v_{2}|$ and $\alpha(u_{1},v_{1}) >
\alpha(u_{2},v_{2})$.
\\
\\
Then there exists $c_{f,g} \in \mathbb{R}^{+}$ such that
\eqref{eq:old317} and \eqref{eq:old318} hold. Specifically,
$c_{f,g}$ is the unique positive root of the equation \be
g_{u_{1},v_{1}} \left( {c^{2} \over u_{1}} \right) = g_{u_{2},v_{2}}
\left( {c^{2} \over u_{2}} \right), \ee where $g_{u,v}(x) :
\mathbb{R}^{+} \rightarrow (0,u+|v|)$ is given by \eqref{eq:old32}.
\end{thm}

\setcounter{equation}{0}
\section{Open Problems}\label{sec:openproblems}

Here is a sample of issues which could be the subject of further
investigations :
\\
\\
\textbf{1.} One unresolved matter is the comparison of arbitrary
difference forms in the range where $N^{-3/4} = O(p)$ and $p =
O(N^{-3/5})$. Here the problem is that the binomial model itself
does not prove of any use. This provides, more generally, motivation
for looking at other models. Obviously one could look at the
so-called \emph{uniform} model on subsets (see \cite{JLR}), but this
seems a more awkward model to handle. Note that the property of one
binary form dominating another is not monotone, or even convex.
\\
\\
\textbf{2.} Secondly, a very tantalizing problem is to investigate
what happens while crossing a sharp threshold, whenever it arises
under the conditions of Theorem \ref{thm:compareforms}(ii).
\\
\\
\textbf{3.} Thirdly, one can ask if the various concentration
estimates in Theorem \ref{thm:mainuniform} can be improved. When $p
= o(N^{-1/2})$ we have only used an ordinary second moment argument,
and it is possible to provide explicit estimates. Explicitly, the
following follows from Chebyshev's Theorem (see the appendix for a proof).

\begin{thm}\label{thm:boundsfnNinpaper} Let $p(N) := cN^{-\delta}$
for some $c > 0$, $\delta \in (1/2,1)$. Set $C := \max(1,c)$, $f(\delta) :=
\min \{ \frac{1}{2}, \frac{3\delta - 1}{2} \}$ and let $g(\delta)$
be any function such that $0 < g(\delta) < f(\delta)$ for all
$\delta \in (1/2,1)$. Set $P_1(N) := (4/c)N^{-(1-\delta)}$ and
$P_2(N) := N^{-(f(\delta)-g(\delta))}$. For any subset chosen with
respect to the binomial model with parameter $p=p(N)$, with
probability at least $1-P_1(N)-P_2(N)$ the ratio of the cardinality
of its difference set to the cardinality of its sumset is $2 +
O_C(N^{-g(\delta)})$. Thus the probability a subset chosen with
respect to the binomial model is not difference dominated is at most
$P_1(N)+P_2(N)$, which tends to zero rapidly with $N$ for $\delta
\in (1/2, 1)$.
\end{thm}

The range $N^{-1/2} = o(p(N))$ seems more interesting, however. Here
we proved that the random variable $\mathscr{S}^{c}$ has expectation
of order $P(N)^{2}$, where $P(N) = 1/p(N)$, and is concentrated
within $P(N)^{3/2} \log^{2} P(N)$ of its mean. Now one can ask whether
the constant $3/2$ can be improved, or at the very least can one get
rid of the logarithm?
\\
\\
\textbf{4.} Finally, it is natural to ask for extensions of our results to
$\mathbb{Z}$-linear forms in more than two variables. Let
\be
f(x_1,...,x_k)\ =\ u_1 x_1 + \cdots + u_k x_k, \;\;\; u_i \in \mathbb{Z}_{\neq 0},
\ee
be such a form. We conjecture the following generalization of Theorem 3.1 :

\begin{conj}\label{conj:linearforms}

Let $p:\N \to (0,1)$ be a function satisfying \eqref{eq:old13}. For a positive
integer $N$, let $A$ be a random subset of $I_{N}$ obtained by choosing each
$n \in I_{N}$ independently with probability $p(N)$. Let $f$ be as in (4.1)
and assume that GCD$(u_1,...,u_n) = 1$. Set
\be
\theta_{f}\ :=\ \# \{\sigma \in S_k : (u_{\sigma(1)},...,u_{\sigma(k)}) =
(u_1,...,u_k) \}.
\ee
Let
$\mathscr{D}_{f}$ denote the random variable $|f(A)|$. Then the
following three situations arise :
\\
\\
(i) $p(N) = o(N^{-1/k})$ : Then \be \mathscr{D}_{f}\ \sim\
{1 \over \theta_{f}} (N \cdot
p(N))^{k}. \ee (ii) $p(N) = c\cdot N^{-1/k}$ for some $c \in
(0,\infty)$ : There is a rational function $R(x_0,...,x_k)$ in $k+1$
variables, which is increasing in $x_0$, and an increasing
function \\ $g_{u_1,...,u_k} :
(0,\infty) \rightarrow (0,\sum_{i=1}^{k} |u_i|)$ such that
\be\label{eq:old44} \mathscr{D}_{f}\ \sim\ g_{u_1,...,u_k} \left(
R(c,u_1,...,u_k) \right) \cdot N. \ee (iii) $N^{-1/k} = o(p(N))$ : Let
$\mathscr{D}_{f}^{c} := \left( \sum_{i=1}^{k} |u_i| \right)N -
\mathscr{D}_{f}$. Then \be
\mathscr{D}_{f}^{c}\ \sim\ {2\theta_{f} \prod_{i=1}^{k} |u_i|
\over p(N)^{k}}. \ee
\end{conj}



\section*{Acknowledgement}

We thank Devdatt Dubhashi for pointing us in the direction of Vu's paper, the
participants of CANT 2007 for interesting conversations, and the referees for
comments on an earlier draft. The second named author was partly supported by NSF
grant DMS0600848.

\appendix

\setcounter{equation}{0}
\section{Explicit bounds}

Here we prove Theorem \ref{thm:boundsfnNinpaper}. The proof uses essentially
only Chebyshev's inequality.
We have deliberately made this section
entirely self-contained, rather than appealing to results from
\cite{AS}, for the benefit of readers who may not be too familiar
with discrete probability theory.

We first establish some notation, and then prove a sequence of
lemmas from which Theorem \ref{thm:boundsfnNinpaper} immediately follows.
Our goal is to provide explicit bounds which decay like $N$ to a
power.

Let $X_{n;N}$ denote the binary indicator variable for $n$ being in
a subset (it is thus $1$ with probability $c N^{-\delta}$ and $0$
otherwise), and let $X$ be the random variable denoting the
cardinality of a subset (thus $X = \sum_n X_{n;N}$). For two pairs
of ordered elements $(m,n)$ and $(m',n')$ in $I_N\times I_N$ ($m <
n$, $m'<n'$), let $Y_{m,n,m',n'} = 1$ if $n-m = n'-m'$, and $0$
otherwise.

\begin{lem}\label{lem:sizeXunif}
With probability at least $1- P_1(N)$, \be X\ \in \ \left[\foh\ c
N^{1-\delta},\ \ \frac32\ c N^{1-\delta} \right]. \ee Let
$\mathcal{O}$ denotes the number of ordered pairs $(m,n)$ (with $m <
n$) in a subset of $I_N$ chosen with respect to the binomial model.
Then with probability at least $1-P_1(N)$ we have
\be\label{eq:numbpairs} \frac{\foh c N^{1-\delta} \left(\foh c
N^{1-\delta} - 1\right)}{2} \ \le \mathcal{O} \ \le \ \frac{\frac32
c N^{1-\delta} \left(\frac32 c N^{1-\delta} - 1\right)}{2}. \ee
\end{lem}

\begin{proof} We have $\E[X] = \sum_n \E[X_{n;N}] = c N^{1-\delta}$. As the $X_{n;N}$ are
independent, \be \sigma^2_{X}\ =\ \sum_n \sigma^2_{X_{n;N}}\ =\ N
\left(cN^{-\delta} - c^2N^{-2\delta}\right). \ee Thus \be \sigma_X \
\le \ \sqrt{c}\cdot N^{\frac{1-\delta}{2}}. \ee By Chebyshev's
inequality, \be {\rm Prob}(|X - cN^{1-\delta}| \le k \sigma_X) \ \ge
\ 1 - \frac1{k^2}. \ee For $X\ \in \ \left[\foh\ c N^{1-\delta},\ \
\frac32\ c N^{1-\delta} \right]$ we choose $k$ so that \be k
\sigma_X \ = \ \foh\ c N^{1-\delta} \ \le \ k
\sqrt{c}N^{\frac{1-\delta}{2}}. \ee Thus $k \ge \foh\ \sqrt{c}
N^{(1-\delta)/2}$, and the probability that $X$ lies in the stated
interval is at least $1 - (c N^{1-\delta}/4)^{-1}$. The second claim
follows from the fact that there are $\ncr{r}{2}$ ways to choose two
distinct objects from $r$ objects.
\end{proof}

\begin{rek} By using the Central Limit Theorem instead of
Chebyshev's inequality we may obtain a better estimate on the
probability of $X$ lying in the desired interval; however, as the
Central Limit Theorem is not available for some of the later
arguments, there is negligible gain in using it here. \end{rek}


\emph{Proof of Theorem 4.1.} By Lemma \ref{lem:sizeXunif},
\eqref{eq:numbpairs} holds with probability at least $1 - P_1(N)$.
The main contribution to the cardinalities of the sumset and the
difference set is from ordered pairs $(m,n)$ with $m < n$. With
probability at least $1 - P_1(N)$ there are on the order
$N^{2-2\delta}$ such pairs, which is much larger than the order
$N^{1-\delta}$ pairs with $m=n$. The proof is completed by showing
that almost all of the ordered pairs yield distinct sums (and
differences). Explicitly, we shall show that for a subset chosen
from $I_N$ with respect to the binomial model, if $\mathcal{O}$ is
the number of ordered pairs (which is of size $N^{2-2\delta}$ with
high probability), then with high probability the cardinality of its
difference set is $2\mathcal{O} + O_C(N^{3-4\delta})$ while the
cardinality of its sumset is $\mathcal{O} + O_C(N^{3-4\delta})$.
This argument crucially uses $\delta > 1/2$ (if $\delta = 1/2$) then
the error term is the same size as the main term, and the more
delicate argument given in the main text is needed). We shall show
that almost all of the ordered pairs generate distinct differences;
the argument for the sums follows similarly.

Each ordered pair $(m,n)$ yields two differences ($m-n$ and $n-m$).
The problem is that two different ordered pairs could generate the
same differences. To calculate the size of the difference set, we
need to control how often two different pairs give the same
differences. Consider two distinct ordered pairs $(m,n)$ and
$(m',n')$ with $m < n$ and $m' < n'$ (as the $N^{1-\delta} \ll
N^{2-2\delta}$ `diagonal' pairs $(n,n)$ yield the same difference,
namely 0, it suffices to study the case of ordered pairs with
distinct elements). Without loss of generality we may assume $m \le
m'$. If $n-m =n'-m'$ then these two pairs contribute the same
differences. There are two possibilities: (1) all four indices are
distinct; (2) $n=m'$.

We calculate the expected number of pairs of non-diagonal ordered
pairs with the same difference by using our binary indicator random
variables $Y_{m,n,m',n'}$. Set \be Y \ = \ \sum_{1 \le m \le m' \le
N}\ \sum_{m' < n' \le N}\ \sum_{m < n \le N \atop n'-m'=n-m}\
Y_{m,n,m',n'}.\ee If the four indices are distinct then
$\E[Y_{m,n,m',n'}] = c^4 N^{-4\delta}$; if $n=m'$ then
$\E[Y_{m,n,m',n'}] = c^3 N^{-3\delta}$.

The number of tuples $(m,n,m',n')$ of distinct integers satisfying
our conditions is bounded by $N^3$ (once $m$, $n$ and $m'$ are
chosen there is at most one choice for $n'\in \{m+1,\dots,N\}$ with
$n'-m'=n-m$)\footnote{Although we do not need the actual value,
simple algebra yields the number of tuples is $N^{3}/6 + O(N^2)$.}.
If instead $n=m'$ then there are at most $N^2$ tuples satisfying our
conditions (once $m$ and $n$ are chosen, $m'$ and $n'$ are uniquely
determined, though they may not satisfy our conditions). Therefore
\be \E[Y] \ \le \ N^3 \cdot c^4 N^{-4\delta} + N^2 \cdot c^2
N^{-3\delta}\ \le \ 2 C^4 N^{3-4\delta} \ee as $\delta \in (1/2,
1)$.

As $N^{3-4\delta}$ is much smaller than $N^{2-2\delta}$ for $\delta
> 1/2$, most of the differences are distinct. To complete the proof,
we need some control on the variance of $Y$. In Lemma
\ref{lem:stddevYum} we show that \be \sigma_Y \ \le \ 7C^4
N^{r(\delta)},\ee where \be 2r(\delta) = \max \{3-4\delta, 5-7\delta
\}. \ee While we cannot use the Central Limit Theorem (as the
$Y_{m,n,m',n'}$ are not independent and also depend on $N$), we may
use Chebyshev's inequality to bound the probability that $Y$ is
close to its mean (recall the mean is at most $2 C^4
N^{3-4\delta}$). We have \be {\rm Prob}(|Y - \E[Y]| \le k \sigma_Y)
\ge 1 - \frac1{k^2}. \ee Simple algebra shows that if we take $k =
N^{2-2\delta-r(\delta)-g(\delta)}$ then with probability at least $1
- N^{-(f(\delta)-g(\delta))}$ we have $Y \le 9C^4
N^{2-2\delta-g(\delta)}$, which is a positive power of $N$ less than
$N^{2-2\delta}$. Thus an at most negligible amount of the
differences are repeated.

The argument for two ordered pairs yielding the same sum proceeds
similarly: if $\mu + \nu = \mu'+\nu'$ then $\nu-\mu' = \nu'-\mu$.

For our ratio to be $2 + O_C(N^{-g(\delta)})$, two events must
happen. As the probability the first does not occur is at most
$P_1(N)$ and the probability the second does not occur is at most
$P_2(N)$, the probability that the two desired events happen is at
least $1-P_1(N)-P_2(N)$.
\\
\\
Except for the claimed estimate on $\sigma_Y$, the above completes
the proof of Theorem \ref{thm:boundsfnNinpaper}. We now prove our bound for
$\sigma_Y$.

\begin{lem}\label{lem:stddevYum} Let the notation be as in
Theorem \ref{thm:boundsfnNinpaper} and (A.10). We have \be \sigma_Y \ \le \
7C^4 N^{r(\delta)}.\ee
\end{lem}

\begin{proof} If $U$ and $V$ are two random variables, then \be {\rm
Var}(U+V) \ = \ {\rm Var}(U) + {\rm Var}(V) + 2{\rm CoVar}(U,V). \ee
By the Cauchy-Schwartz inequality, ${\rm CoVar}(U,V) \le \sqrt{{\rm
Var}(U){\rm Var}(V)}$. Thus \be {\rm Var}(U+V) \ \le \ 3 {\rm
Var}(U) + 3{\rm Var}(V). \ee We may therefore write \be \sum
Y_{m,n,m',n'} \ = \ \sum U_{m,n,m',n'} + \sum V_{m,n,n'} \ = \ U +
V, \ee where in the $U$-sum all four indices are distinct (with $1
\le m < m' \le N$, $m < n \le N$, $m' < n' \le N$ and $n-m = n'-m'$)
and in the $V$-sum all three indices are distinct (with $1 \le m < n
< n' \le N$ and and $n-m = n'-n$). As ${\rm Var}(Y) \le 3 {\rm
Var}(U) + 3 {\rm Var}(V)$, we are reduced to bounding the variances
of $U$ and $V$.

We first bound ${\rm Var}(U)$. Standard algebra yields \bea {\rm
Var}(U)& \ = \ & {\rm Var}\left(\sum U_{m,n,m',n'}\right)\nonumber\\
& \ =\ & \sum {\rm Var}(U_{m,n,m',n'}) + 2 \sum_{(m,n,m',n') \neq
(\widetilde{m},\widetilde{n},\widetilde{m}',\widetilde{n}')} {\rm
CoVar}(U_{m,n,m',n'},U_{\widetilde{m},\widetilde{n},\widetilde{m}',\widetilde{n}'}).
\nonumber\\ \eea As ${\rm Var}(U_{m,n,m',n'}) = c^4 N^{-4\delta} -
c^8 N^{-8\delta}$ and there are at most $N^3$ ordered tuples
$(m,n,m',n')$ of distinct integers with $n-m=m'-n'$, the ${\rm
Var}(U_{m,n,m',n'})$ term is bounded by $c^4 N^{3-4\delta}$.

For the covariance piece, if all eight indices
($m,n,m',n',\widetilde{m},\widetilde{n},\widetilde{m}',\widetilde{n}'$)
are distinct, then $U_{m,n,m',n'}$ and
$U_{\widetilde{m},\widetilde{n},\widetilde{m}',\widetilde{n}'}$ are
independent and thus the covariance is zero. There are four cases;
in each case there are always at most $N^3$ choices for the tuple
$(m,n,m',n')$, but often there will be significantly fewer choices
for the tuple
$(\widetilde{m},\widetilde{n},\widetilde{m}',\widetilde{n}')$. We
only provide complete details for the first and third cases, as the
other cases follow similarly.

\bi
\item Seven distinct indices: There are at most $N^2$ choices for
$(\widetilde{m},\widetilde{n},\widetilde{m}',\widetilde{n}')$. The
covariance of each such term is bounded by $c^7 N^{-7\delta}$. To
see this, note \bea & & {\rm
CoVar}(U_{m,n,m',n'},U_{\widetilde{m},\widetilde{n},\widetilde{m}',\widetilde{n}'})
\nonumber\\ & & \ \ \ \ \ \ =\
\E[U_{m,n,m',n'}U_{\widetilde{m},\widetilde{n},\widetilde{m}',\widetilde{n}'}]
-
\E[U_{m,n,m',n'}]\E[U_{\widetilde{m},\widetilde{n},\widetilde{m}',\widetilde{n}'}].\eea
The product of the expected values is $c^8 N^{-8\delta}$, while the
expected value of the product is $c^7 N^{-7\delta}$. Thus the
covariances of these terms contribute at most $c^7 N^{5-7\delta}$.

\item Six distinct indices: The
covariances of these terms contribute at most $c^6 N^{4-6\delta}$.

\item Five distinct indices: The
covariances of these terms contribute at most $c^5 N^{3-5\delta}$
(once three of the
$\widetilde{m},\widetilde{n},\widetilde{m}',\widetilde{n}'$  have
been determined, the fourth is uniquely determined; thus there are
at most $N^3$ choices for the first tuple and at most 1 choice for
the second).

\item Four distinct indices: The
covariances of these terms contribute at most $c^4 N^{3-4\delta}$.

\ei The $N$-dependence from the case of seven distinct indices is
greater than the $N$-dependence of the other cases (except for the
case of four distinct indices if $\delta
> 2/3$). We also only increase the contributions if we replace $c$ with
$C = \max(c,1)$. We therefore find \bea {\rm Var}(U) & \ \le \ & C^4
N^{3-4\delta} + 2\left(C^7 N^{5-7\delta} + C^6N^{4-6\delta} + C^5
N^{3-5\delta}+ C^4 N^{3-4\delta}\right) \nonumber\\ & = & 3C^4
N^{3-4\delta} + 6C^7 N^{5-7\delta}. \eea

Similarly we have \bea {\rm Var}(V)& \ = \ & {\rm Var}(\sum
V_{m,n,n'})\nonumber\\ & \ =\ & \sum {\rm Var}(V_{m,n,n'}) + 2
\sum_{(m,n,n') \neq (\widetilde{m},\widetilde{n},\widetilde{n}')}
{\rm
CoVar}(V_{m,n,n'},V_{\widetilde{m},\widetilde{n},\widetilde{n}'}).
\nonumber\\ \eea The ${\rm Var}(V_{m,n,n'})$ piece is bounded by
$N^2 \cdot c^3 N^{-3\delta}$ (as there are at most $N^2$ tuples with
$n'-n = n-m$). The covariance terms vanish if the six indices are
distinct. A similar argument as before yields bounds of
$c^5N^{3-5\delta}$ for five distinct indices, $c^4 N^{2-4\delta}$
for four distinct indices, and $c^3 N^{2-3\delta}$ for three
distinct indices. The largest $N$-dependence is from the $c^3
N^{2-3\delta}$ term (as $\delta > 1/2$). Arguing as before and
replacing $c$ with $C$ yields \bea {\rm Var}(V) \ \le \ C^3
N^{2-3\delta} + 2\cdot 3C^3 N^{2-3\delta} \ \le \ 7C^3
N^{2-3\delta}. \eea As $\delta < 1$, $2-3\delta < 3 - 4\delta$.
Therefore \bea {\rm Var}(Y) & \ \le \ & 3\cdot \left(3C^4
N^{3-4\delta} + 6C^7 N^{5-7\delta}\right) + 3 \cdot
7C^3 N^{2-3\delta} \nonumber\\
& \le & 30C^4 N^{3-4\delta} + 18C^7 N^{5-7\delta} \ \le \ 49C^8
N^{2r(\delta)}, \eea which yields \be \sigma_Y \ \le \ 7C^4
N^{r(\delta)}.\ee
\end{proof}

\begin{rek}
An extreme choice of $g$ would be to choose $g(\delta) = \epsilon$,
for some small positive constant $\epsilon$. Since $f(\delta) \geq
1/4$ for all $\delta \in (1/2,1)$, we then obtain a bound of
$2+O_{C}(N^{-\epsilon})$ for the ratio of the cardinality of the
difference set to the sumset with probability $1 -
O_{C}(N^{-\min\{1-\delta, \frac{1}{4}-\epsilon\}})$.
\end{rek}

\begin{rek}
Alternatively, one can get a tighter bound on the ratio than in
Theorem \ref{thm:boundsfnNinpaper} at the expense of having a bound on the
probability which is non-trivial only for $\delta < \delta_{\max} <
1$. For example, if one instead chooses $k =
N^{3-4\delta-r(\delta)}$ in (A.11) and $P_{2}(N) :=
N^{-(6-8\delta-2r(\delta))}$, then the statement of Theorem
\ref{thm:boundsfnNinpaper} still holds, but is a non-trivial statement only
for $6-8\delta-2r(\delta) > 0$, i.e.: for $\delta < 3/4$. This is a
natural choice of $\delta_{\max}$, given the results of \cite{GJLR}.
Specifically, if $\delta > 3/4$ then almost surely no differences or
sums are repeated, and the set is a Sidon set (and therefore we
trivially have the ratio of the cardinality of the difference set to
the sumset is approximately $2$).
\end{rek}


\ \\

\end{document}